\documentclass[12pt,a4paper,reqno]{amsart}

\usepackage{amsfonts,amsmath}
\usepackage{graphicx}
\usepackage{anysize}
\usepackage[colorlinks,citecolor={red}]{hyperref}

\usepackage{tikz,xcolor}

\definecolor{lime}{HTML}{A6CE39}
\DeclareRobustCommand{\orcidicon}
	{
		\begin{tikzpicture}
		\draw[lime, fill=lime] (0,0) 
		circle [radius=0.17] 
		node[white] {{\fontfamily{qag}\selectfont \tiny ID}};
		\draw[white, fill=white] (-0.0625,0.095) 
		circle [radius=0.007];
		\end{tikzpicture}
		\hspace{-2mm}
	}
\foreach \x in {A, ..., Z}
	{\expandafter\xdef\csname orcid\x\endcsname
		{\noexpand\href{https://orcid.org/\csname orcidauthor\x\endcsname}
			{\noexpand\orcidicon}
		}
	}

\marginsize{3cm}{3cm}{3cm}{3cm}
\newtheorem{theorem}{Theorem}[section]
\newtheorem{remark}{Remark}[section]
\newtheorem{case}{Case}[section]

\newtheorem{example}{Example}[section]

\begin{document}
\title{A new technique to solve linear integro-differential equations (IDEs) with modified Bernoulli polynomials}
\author[Udaya Pratap Singh] {Udaya Pratap Singh\hspace{-0.5em}\orcidA{} \\ \tiny{Department of Applied Sciences, Rajkiya Engineering College,\\ Sonbhadra, Uttar Pradesh, India}}
%
\footnote{email: \href{}{upsingh1980@gmail.com}\\ Orcid\hspace*{-0.5em}\orcidA{}: \url{https://orcid.org/0000-0002-4538-9377} }
	\begin{abstract}
		 In this work, a new technique has been presented to find approximate solution of linear integro-differential equations. The method is based on modified orthonormal Bernoulli polynomials and an operational matrix thereof. The method converts a given integro-differential equation into a set of algebraic equations with unknown coefficients, which is easily obtained with help of the known functions appearing in the equation, modified Bernoulli polynomials and operational matrix. Approximate solution is obtained in form of a polynomial of required degree. The method is also applied to three well known integro-differential equations to demonstrated the accuracy and efficacy of the method. Numerical results of approximate solution are plotted to compared with available exact solutions. Considerably small error of approximation is observed through numerical comparison, which is further reducible to a required level of significance. Method is comparatively simpler and shorter than many existing methods.
		
	\end{abstract}
	\maketitle 
	\bigskip

	\section{Introduction}
	\label{section 1}
	
	Mathematical modeling of real world problems often give rise to ordinary or partial differential equations, integral equations, ordinary or partial integro-differential equations and some other forms. Of these, integro-differential equations (IDEs) appear in almost all areas of science and engineering. Mathematical formulations of physical phenomena such as population problem, nono-hydrodynamics, fluid mechanics biological and chemical models, ecology models, financial problem, process engineering, aerospace and design engineering, hydro-electric machines, reactor dynamics, and many more are examples where IDEs are frequently encountered. 
	
	Many of such IDEs are difficult to solve for analytic solution and require an efficient approximation or numerical technique. Solution to IDEs in different fields have been point of attracting attention of researchers from long past \cite{McKendrick1926, Schaffer1986} including some notable contributions on mathematical modeling of the spread of infection \cite{Kendall1965}, problems in quantum mechanics \cite{Baym1990}, detailed consideration of integro-differential equations theory \cite{Vangipuram1995}, problems of hydrodynamics with incompressible viscous fluids \cite{Petrin1997} and mathematical formulation in ecology \cite{Kot2001}, however, researchers have focused towards numerical techniques to solve such problems with the evolution of computers science. 
	
	Bulk of recent literature is available to explore approximate and numerical solutions of IDEs \cite{Darania2007,Abdelaziz2010,Yuzbasi2011,Sweilam2013}. Some latest contributions on numerical solution of IDEs in recent years are meshless method \cite{Dehghan2012}, Bernstein operational matrix approach \cite{Maleknejad2012}, collocation approach \cite{Yuzbasi2012}, improved Legendre method \cite{Yuzbasi2013}, with Euler polynomials \cite{Mirzaee2015}, operational matrices method \cite{Singh2017b}, Taylor collocation approach \cite{Jafarzadeh2018}, convolution integrals approach \cite{Katsikadelis2019} and that involving special functions \cite{Singh2020}.
	
	Numerical and approximation techniques for ODEs, PDEs and IDEs also involved different well known polynomials, such as applications of Bernoulli polynomials \cite{Cheon2003}, Chebyshev polynomials approach \cite{Maleknejad2007}, application of Legendre polynomials \cite{Nemati2015}, Laguerre polynomials and Wavelet Galerkin method \cite{Rahman2012}, Legendre wavelets \cite{Yousefi2006}, the operational matrix \cite{Sahu2019}. 

	Amongst these, many authors also used Bernoulli polynomials in different ways to find numerical solution of many complex problems, such as numerical approximation for generalized pantograph equation using Bernoulli matrix method \cite{Tohidi2013}, numerical solution of second-order linear system of partial differential equations using Bernoulli polynomials \cite{Tohidi2013a}, numerical solution of Volterra type integral equations by means of Bernoulli polynomials \cite{Mohsenyzadeh2016}. 
	
	In this work, it is proposed to find polynomial approximation to solution of linear integro-differential equations (IDEs) by application of an operational matrix developed from a class of modified Bernoulli polynomials.
	
	\section{Bernoulli Polynomials}
	\label{section 2: Bernoulli Polynomials }
	Jakob Bernoulli, in late seventeenth century, discussed some special polynomials in his book "$Ars Conjectandi$", which were explicitly studied by Leonhard Euler, who established the finite difference relation ${B_m}(x  + 1) - {B_m}(x ) = m{x ^{m - 1}},\,m \ge 1$ for these polynomials in his book “Foundations of differential calculus” in 1755, and also suggested the method of generating function to derive these polynomials \cite{Costabile2006}. Later on, J. L. Raabe in 1851 discussed these polynomials together with Bernoulli numbers in connection with formula $\sum\limits_{n = 0}^{m - 1} {{B_n}\left( {x + \frac{k}{m}} \right)}  = {m^{ - (n + 1)}}{B_n}(mx)$ and termed the polynomials $B_n(x)$ as $Bernoulli \, Polynomial$. After Raabe, many researchers paid their attention towards properties of these polynomials. The most common formula for Bernoulli polynomials in recent mathematical applications is: 
	\begin{equation}\label{eq.1 : Basic Bernoulli Polynomials}
	{B_n}(\zeta ) = \sum\limits_{j = 0}^n {\,\left( {\begin{array}{*{20}{c}}
			n \\ j 	\end{array}} \right)\,\,{B_j}(0)\,{\zeta ^{n - j}},\hspace{1em} n = 0,1,2,...} \hspace{1em} 0 \le \zeta  \le 1
	\end{equation}
	which was presented by Costabile and Dell’Accio \cite{Costabile2006}. The numbers $B_j(0)$ are the Bernoulli numbers, which can also be calculated with Kronecker’s formula \cite{Todorov1984}:
	 
	\begin{equation}
	{B_n}(0) =  - \sum\limits_{j = 1}^{n + 1} {\frac{{{{( - 1)}^j}}}{j}\left( \begin{array}{c}
		n + 1\\	j
		\end{array} \right)} \sum\limits_{k = 0}^j {{k^n}} ; \: n \geq 0
	\end{equation}
	
	For illustration, first few Bernoulli polynomials can be written as ${B_0}(x) = 1, B_1(x)=x-\frac{1}{2}, B_2(x)=x^2-x+\frac{1}{6}	, 	B_3(x)=x^3-\frac{3}{2}x^2+\frac{1}{2}x	, 	B_4(x)=x^4-2x^3+x^2-\frac{1}{30}$.
	
	An interesting property of Bernoulli polynomials is that they form a complete basis over $[0,1]$ \cite{Kreyszig1978}, which extends the applicability of these polynomials towards various numerical approximation techniques. In the present work, this property will be used as an underlying characteristics together with the following key relations of Bernoulli polynomials \cite{Costabile2006} :
	\begin{equation}\label{eq.4 : properties of Bernoulli Polnmls}
	\left. \begin{gathered}
	B'_n(\zeta ) = n{B_{n - 1}}(\zeta ),\,\,\,n \ge 1 \hfill \\
	\int_0^1 {{B_n}(z)dz = 0,\,\,\,\,\,\,\,\,\,\,n \ge 1}  \hfill \\
	\end{gathered}  \right\}.
	\end{equation}
	
	Some other properties, generalization and applications of Bernoulli polynomials can be found in notable literature  \cite{Kurt2011, Natalini2003, Lu2011}.
	
	\section{Modified Bernoulli Polynomials}
	\label{section 3: The Orthonormal Polynomials}
	It can be easily verified that the polynomials ${B_n}(x)\,(n \ge 1)$  given by equation (\ref{eq.1 : Basic Bernoulli Polynomials})  are orthogonal to $B_o(x)$ with respect to standard inner product on  ${L^2} [0,1]$ :
	\begin{equation}
	<f_1,f_2> = \int_{0}^{1} f_1(x) \bar{f_2(x)} dx \,; \:\: f_1,f_2 \in L^2[0,1]
	\end{equation}
	
	 Using this property, an orthonormal set of polynomials can be derived for any $B_n(x)$ with Gram-Schmidt orthogonalization. First few of such orthonormal polynomials are obtained as:
	\begin{subequations} \label{Orthonormal polynomials}
	\begin{equation} \label{Orthonormal degree 0}
	{\phi _{0\,}}(x ) = 1
	\end{equation}
	\begin{equation}
	{\phi _1}(x ) = \sqrt 3 ( - 1 + 2x )
	\end{equation}
	\begin{equation}
	\phi_2\left(x\right)=\sqrt5\left(1-6x+6x^2\right)
	\end{equation}
	\begin{equation}
	\phi_3(x)=\sqrt7(-1+12x-30x^2+20x^3)
	\end{equation}
	\begin{equation}
	\phi_4(x)=3(1-20x+90x^2-140x^3+70xa^4)
	\end{equation}
	\begin{equation}
	\phi_5(x)=\sqrt{11}(-1+30x-210x^2+560x^3-630x^4+252 x^5)
	\end{equation}
	\begin{equation}
	{\phi _6}(x) = \sqrt {13}\, ( 1 - 42x  + 420{x ^2} - 1680{x ^3}
	+ 3150{x ^4} - 2772 x^5 + 924{x ^6})
	\end{equation}
	\begin{equation}
	\begin{split}
	{\phi _7}\left( x \right) = \sqrt {15}\, (- 1 + 56x &- 756{x^2} + 4200{x^3}\\
	&- 11550{x^4} + 16632{x^5} - 12012{x^6} + 3432{x^7})
	\end{split}
	\end{equation}

	\end{subequations}
	\section{Approximation of Functions}
	\label{section 4 : Approximation of Functions}
	\begin{theorem} \label{theorem 1}
		Let $H=L^2[0,1]$  be a Hilbert space and $Y=span\left\{y_0,y_1,y_2,...,y_n\right\}$  be a subspace of $H$ such that $dim{(}Y)<\infty$ , every $f\in H$  has a unique best approximation out of $Y$ \cite{Kreyszig1978}, that is, $\forall y(t)\in Y,\, \exists \, \hat{f}(t)\in Y$ s.t. $\parallel f(t)-\hat{f}(t)\parallel_2\le\parallel f(t)-y(t)\parallel_2$. This implies that,  $\forall \, y(t)\in Y, <f(t)-f(t), y(t)>= 0$, where $<,>$  is standard inner product on $L^2\in[0,1]$ (\textit{c.f. Theorems 6.1-1 and 6.2-5, Chapter 6} \cite{Kreyszig1978}).
	\end{theorem}
	\begin{remark}
		Let $Y=span\left\{\phi_0,\phi_1,\phi_2,...,\phi_n\right\},$ where $\phi_k\in L^2[0,1]$ are orthonormal Bernoulli polynomials. Then, from Theorem \ref{theorem 1}, for any function $ f\in L^2[0,1],$
		\begin{equation} \label{eq.5 : approx theorem}
		f\approx\hat{f}=\sum_{k=0}^{n}{c_k\phi_k}, 
		\end{equation}
		where $c_k=\left\langle f,\phi_k\right\rangle,$ and $<,>$ is the standard inner product on $L^2\in[0,1]$.
	\end{remark}
	For numerical approximation, series (5) can be writequation as:
	\begin{equation} \label{eq.6 : approximation series}
	f(x)\simeq\sum_{k=0}^{n}{c_k\phi_k(x)=C^T\phi(x)} 
	\end{equation}
	where $C=\left(c_0,c_1,c_2,...,c_n\right), \phi(x)=\left(\phi_0,\phi_1,\phi_2,...,\phi_n\right)$ are column vectors. The number of polynomials $n$ can be chosen to meet required accuracy.
	\section{Construction of operational matrix}
	\label{section 5 : Construction of operational matrix}
	The orthonormal polynomials, as derived in section \ref{section 3: The Orthonormal Polynomials}, can be integrated as follows:
	\begin{equation}\label{eq.7 : Ortho Pol. of deg. 0}
	\int_0^x  {{\phi _o}(t)dt }  = \frac{1}{2}{\phi _o}(x) + \frac{1}{{2\sqrt 3 }}{\phi _1}(x)
	\end{equation}
	\begin{equation}\label{eq.8 : Ortho Pols.}
	\begin{array}{l}
	\int\limits_0^x  {{\phi _i}(t)dt = } \,\,\,\,\,\frac{1}{{2\sqrt {(2i - 1)(2i + 1)} }}{\phi _{i - 1}}(x)\\
	\hspace*{4.7em} + \frac{1}{{2\sqrt {(2i + 1)(2i + 3)} }}{\phi _{i + 1}}(x);\hspace{0.5em} (i = 1\,,2,...\,,n)
	\end{array}
	\end{equation}
	Relations (\ref{eq.7 : Ortho Pol. of deg. 0}-\ref{eq.8 : Ortho Pols.}) are combined to closed form as:
	\begin{equation}\label{Combined Expression for Orthonormal Polynomials}
	\int\limits_0^x  {\phi(\eta )d\eta  = \,\,} \Theta \,\phi(x)
	\end{equation}
	where $x\in[0,1]$ and $\Theta$ is operational matrix of order $(n+1)$ given as :
	\begin{equation} \label{Operational Matrix}
	\Theta \, = \frac{1}{2}\left[ {\begin{array}{*{20}{c}}
		1&{\frac{1}{{\sqrt {1.3} }}}&0& \cdots &0\\
		{ - \frac{1}{{\sqrt {1.3} }}}&0&{\frac{1}{{\sqrt {3.5} }}}& \cdots &0\\
		0&{ - \frac{1}{{\sqrt {3.5} }}}&0& \ddots & \vdots \\
		\vdots & \vdots & \ddots & \ddots &{\frac{1}{{\sqrt {\left( {2n - 1} \right)\left( {2n + 1} \right)} }}}\\
		0&0& \cdots &{ - \frac{1}{{\sqrt {\left( {2n - 1} \right)\left( {2n + 1} \right)} }}}&0
		\end{array}} \right]
	\end{equation}
	\section{Solution of Integro-differential Equation}
	\label{section : Solution of IDE}
	The relation (\ref{Combined Expression for Orthonormal Polynomials}) together with the operational matrix (\ref{Operational Matrix}) prove to be a handy tool to solve linear integro-differential equations of the form:
	\begin{equation} \label{IDE Operator Form}
	\mathbb{L}y(x) + f(x) \int_{0}^{x} K(x,t)\, y^{j}(t)\,dt; \: x\in[0,1]
	\end{equation}
	\noindent with suitable initial conditions on $y$ and its derivatives, where $\mathbb{L}$ is linear differential operator of order $k$, $y^{j}$ denotes $j^{th}$ derivative of $y$ for $j<k$, $f(x)$ is some continuous function of $x$ and $K(x,t)$ is non-singular kernel of integration. 
	
	However, the method works for any integro-differential equation of type (\ref{IDE Operator Form}), it is fast and sophisticated for integro-differential equations of the form:
	\begin{equation}\label{IDE Simple form}
	\mathbb{L}y(x) +\mathbb{\mu} \int_{0}^{x}(x-t)^{m-1}y^{j}(t)dt = r(x)
	\end{equation}
	\noindent where, $\mu$ is constant and $m$ is some finite positive integer.
	
	In order to present basic steps of the method in simpler way, we will first consider a simple integro-differential equation of from (\ref{IDE Simple form}). Solution to general form (\ref{IDE Operator Form}) will be presented subsequently with necessary modifications.
	

	\begin{case} \label{Case 1}
		Linear integro-differential equations with constant coefficients
	\end{case}
	Let us take the second order linear integro-differential equations of the form:
	\begin{equation}\label{IDE const coeff}
	\begin{split}
	a_k\frac{d^ky}{dx^k} +a_{k-1}\frac{d^{k-1}y}{dx^{k-1}} &+\dots + a_2\frac{d^2y}{dx^2} + a_1 \frac{dy}{dx} \\&+  a_o y + b \int_{0}^{x}(x-t)^{m-1}y^j(t)dt = r(x)\\
	y(0) =y_0,&\, \left( \frac{dy}{dx} \right)_{x=0} = y_1,\, \dots, \, \left( \frac{d^{k-1}y}{dx^{k-1}} \right)_{x=0}=y_{k-1} 
	\end{split}
	\end{equation}
	
	\noindent where, $m$ is some finite positive integer and $y^j = \frac{d^jy}{dx^j}$ for $j\le k$. It is assumed that $r(x)$ is continuous and equation (\ref{IDE const coeff}) admits a unique solution on $[0,1]$.
	
	Let us take 
	\begin{equation}\label{y(n)= C^T.PHI}
	\frac{d^ky}{dx^k} = C^T\,\phi(x)
	\end{equation}
	\noindent so that equation (\ref{IDE const coeff}) takes the form:
	\begin{equation}\label{IDE const with yn=C^T phi}
	\begin{split}
	a_k& C^T \phi(x) +a_{k-1}C^T \Theta \phi(x) +\dots+ a_2 C^T \Theta^{k-2} + a_1 C^T \Theta^{k-1} \\&+  a_o C^T \Theta^k\phi(x) + b C^T \Theta^{k-j} \int_{0}^{x}(x-t)^{m-1}\phi(t)dt = p(x)+r(x) = \tilde{r}(x) 
	\end{split}
	\end{equation}
	\noindent where, $p(x)$ is a polynomial of degree $m+k$ arising due to initial conditions. 
	
	Now, noting that $\int_{0}^{x}(x-t)^{m-1}\phi(t)dt = \Theta^m\phi(x)$, and taking
	\begin{equation}\label{r(x)= R^T.PHI}
	\tilde{r}(x) = R^T\,\phi(x)
	\end{equation}
	for some real vector $R^T=(r_o,r_1,\dots,r_n)$ of dimension $1\times(n+1)$, equation (\ref{IDE const with yn=C^T phi}) simplifies to- 
	\begin{equation}\label{IDE Simple Transformed}
	C^T \, \left( a_k\, \mathbb{I} + a_{k-1}\,{\Theta} +\dots + a_1 \, \Theta^{k-1} + a_o\,\Theta ^k +  b\,{\Theta ^{m+k-j}} \right)\,\phi(x ) = R^T\,\phi (x)
	\end{equation}
	
	Equation (\ref{IDE Simple Transformed}) gives:
	\begin{equation}\label{C^T for IDE const coeff.}
	C^T = R^T\,\left( a_k\, \mathbb{I} + a_{k-1}\,{\Theta} +\dots + a_1 \, \Theta^{k-1} + a_o\,\Theta ^k +  b\,{\Theta ^{m+k-j}} \right)^{-1}
	\end{equation}
	
	Equation (\ref{C^T for IDE const coeff.}) with equation (\ref{y(n)= C^T.PHI}) give an approximate solution for Integro-differential equation (\ref{IDE const coeff}) as a polynomial of degree $n$ as follows:
	\begin{equation}\label{soln. y(x) for IDE const coeff.}
	y(x) =  y_0 + y_1 x + \frac{y_2}{2!}x^2+\dots+y_{k-1} \frac{y{k-1}}{(k-1)!}x^{k-1} + {C^T\,\Theta ^{k-j}}\,\phi(x)
	\end{equation}
	\noindent where, $\mathbb{I}$ is identity matrix of dimension $n+1$.
	
	%
	\begin{case} \label{Case 2}
	Linear integro-differential equations  with variable coefficients
	\end{case}
	To extend the method discussed in case \ref{Case 1}	for linear integro-differential equations of type (\ref{IDE Operator Form}) with variable coefficients, let us consider the following second order integro-differential equation.
	\begin{equation}\label{IDE General}
	\begin{split}
	\frac{d^2y}{dx^2} + a_1(x) \frac{dy}{dx} +  &a_0(x) y + f(x) \int_{0}^{x}K(x,t)y^k(t)dt = r(x)\\
	&y(0)=\alpha, \: \left(\frac{dy}{dx}\right)_{x=0} = \beta.
	\end{split}
	\end{equation}
	
	\noindent where, $a_0, a_1, f$ and $r$ are continuous functions of $x$, $K(x,t)$ is non-singular kernel and equation (\ref{IDE General}) admits a unique solution on $[0,1]$, and $y^k(x) = \frac{d^ky}{dx^k}$ for $k= 0, \, 1 \text{ or } 2$.
	
	Taking,
	\begin{equation}\label{y''= C^T.PHI}
	\frac{d^2y}{dx^2} = C^T\,\phi(x)
	\end{equation}
	equation  (\ref{IDE General}) can be written as:
	\begin{equation} \label{Integral Breaking}
	\begin{split}
	C^T \phi(x) &+ a_1(x) C^T \Theta \phi(x) + a_0(x) C^T \Theta^2 \phi(x) \\ &+ C^T \Theta^{2-k}  f(x) \int_{0}^{x}K(x,t)\phi(t)dt = g(x)+r(x) = \tilde{r}(x)
	\end{split}
	\end{equation}
	 \noindent where,
	 \begin{equation}
	 g(x)={a_1}\alpha  + {a_0}\left( {\alpha  + \beta x} \right) + 
	 \left\{ 
	 	{\begin{array}{*{20}{c}}
	 	{f(x) \int_{0}^{x}K(x,t)(\alpha  + \beta t);\:k = 0}\\
	 	{\beta f(x) \int_{0}^{x}K(x,t) ;\:k = 1}\\
	 	{0;\:k = 2}
	 	\end{array}} \right.
	 \end{equation}	 
	The integral terms on left side of equation (\ref{Integral Breaking}) can be approximated by modified Bernoulli polynomials as $f(x) \int_{0}^{x}K(x,t)\phi(t)dt = \hat{\Theta}\phi(x)$  for some matrix $\hat{\Theta}$ of order $(n+1)$. Now, writing $\Theta^{2-k} \hat{\Theta} = \Phi $ and taking $\tilde{r}(x) = R^T\phi(x)$, equation  (\ref{Integral Breaking}) can be written as 
	\begin{equation}\label{IDE General Transformed Re-arranged}
	\begin{split}
	{C^T}{\phi}(x) +  C^T\Theta \left[ a_1(x) \phi(x) \right] + C^T \Theta^2 \left[ a_0(x) \phi(x)\right] + C^T \Phi \phi(x) = R^T\,{\phi}( x )
	\end{split}
	\end{equation}
	\noindent where, $a_1(x) \phi(x)$ and $a_0(x) \phi(x)$ are column vectors of type 
	\begin{equation} \label{Merging of aT & phi }
	\begin{split}
	a_k(x) \phi(x) &= \left( a_k(x)  \phi_0(x), \,a_k(x)  \phi_1(x), \, \dots , a_k(x)  \phi_n(x) \right) \\
	&\equiv (\psi_{ko}(x),\psi_{k1}(x),\dots,\psi_{kn}(x)) = \psi_k(x); \: k=0,1. \hspace{5mm}(say!)
	\end{split}	
	\end{equation} 
	\noindent of which, each $\psi_{kj}(x)$ can be approximated as a linear combination of orthonormal polynomials in the form $ \psi_{kj}(x) = A_{kj}^T \phi(x)$ such that $A_{kj}^T$ are vectors of form $1\times (n+1)$ for $ k=0,1$ and $j = 1,2,\dots,n$. 
	
	Therefore,  $a_k(x) \phi(x) = \psi_k(x) = \hat{A_k} \phi(x)$, where $\hat{A_k} = \left( A_{k0}^T,A_{k1}^T, \dots, A_{kn}^T \right)$ are matrices of dimension $(n+1)\times(n+1)$ for $k=0,1$. Finally, using these intermediate approximations, equation (\ref{IDE General Transformed Re-arranged}) can be written as:
	\begin{equation}\label{IDE General final form}
	C^T \left( \mathbb{I} + \Theta \hat{A_1}^T + \Theta^2 \hat{A_0}^T + \Phi \right) \phi(x)= R^T\,{\phi}(x)
	\end{equation}
		
	Equation (\ref{IDE General final form}), gives:
	\begin{equation}\label{C^T for Gen Form}
	C^T = R^T \left( \mathbb{I} + \Theta \hat{A_1}^T + \Theta^2 \hat{A_0}^T + \Phi \right)^{-1}
	\end{equation}
	
	\noindent where, $\mathbb{I}$ is identity matrix of order $n+1$. Now, using equation (\ref{C^T for Gen Form}) in equation (\ref{y''= C^T.PHI}), the expression for $y(x)$ is obtained as:
	\begin{equation}
	y(x) = \alpha + \beta x + R^T \left( \mathbb{I} + \Theta \hat{A_1}^T + \Theta^2 \hat{A_0}^T + \Phi \right)^{-1} \,\Theta^2\,{\phi}(x)
	\end{equation}
	%
	\section{Examples}
	\label{section 7 : Examples}
	In order to establish the accuracy and efficacy of the method discussed so far, examples have been taken from earlier proved results. 
	\begin{example} \label{example 1}
		Let us consider the integro-differential equation of fourth order \cite{Sweilam2013}: 
		\begin{equation} \label{Ex.1}
		\begin{split}
			\frac{d^4y}{dx^4} - y + \int_{0}^{x}y(t)dt =& x + (x+3)e^x\\
		y(0)=1, \: y'(0)=1,  & y''(0)=2, \: y'''(0)=3; \: x \in [0,1]
		\end{split}
		\end{equation}
		
		which has exact solution $y(x) = 1 + x \, e^x$.
	\end{example}

	Let us assume that the solution $y(x)$ can be approximated by modified Bernoulli polynomials of degree $0$ through $7$ such that$-$ 
	\begin{equation}
	\frac{d^4y}{dx^4} = C^T \phi(x)
	\end{equation}
	\noindent where, $\phi(x) = (\phi_0,\phi_1,\phi_2,\phi_3,\phi_4,\phi_5,\phi_6,\phi_7)(x)$ and $C = (c_0,c_1,c_2,c_3,c_4,c_5.c_6,c_7)$ is unknown vector to be determined. Then, equation (\ref{Ex.1}) can be written as$-$ 
	\begin{equation}\label{simplified LHS ex.1}
	C^T\phi(x) - C^T \Theta^4 \phi(x) + C^T \Theta^5 \phi(x) = 1 + x + \frac{x^2}{2} + \frac{x^3}{6} - \frac{x^4}{8} + (x+3) e^x
	\end{equation}
	
	Approximating right side of equation (\ref{simplified LHS ex.1}) as $R^T \phi(x)$, we get 
	\begin{equation}\label{simplified RHS ex.1}
	C^T\phi(x) - C^T \Theta^4 \phi(x) + C^T \Theta^5 \phi(x) = R^T\phi(x)
	\end{equation}
	\noindent where, 
	\begin{equation}
	\begin{split}
	R^T=(7.83814, 2.6674, 0.386136, 0.0327736, 0.00188342, 0.000138055,\\	5.83649 \times 10^{-6}, 1.49271 \times 10^{-7}).
	\end{split} 
	\end{equation}
	
	The unknown vector $C$ and approximate solution $y(x)$ are obtained as:
		\begin{equation}
		\begin{split}
		C^T =& R^T(\mathbb{I} -\Theta^4 + \Theta^5)^{-1}\\ =&(7.87309,\: 2.7079,\:  0.408755,\:  0.039614,\:  0.00280952,\:  \\ &0.000154652,\: 6.52936 \times 10^{-6},\:  1.74013 \times 10^{-7})
		\end{split} 
		\end{equation}
		\begin{equation} \label{y(x) for Ex.1}
		\begin{split}
		y(x) \approx 1 + x + 0.999967 x^2 + 0.500239 x^3 + 0.16579 x^4 \\ + 0.0434333 x^5 + 0.00637214 x^6 + 0.00247836 x^7
		\end{split}
		\end{equation}
		
\begin{figure}[h]
			\centering
			\includegraphics[width=0.49\linewidth]{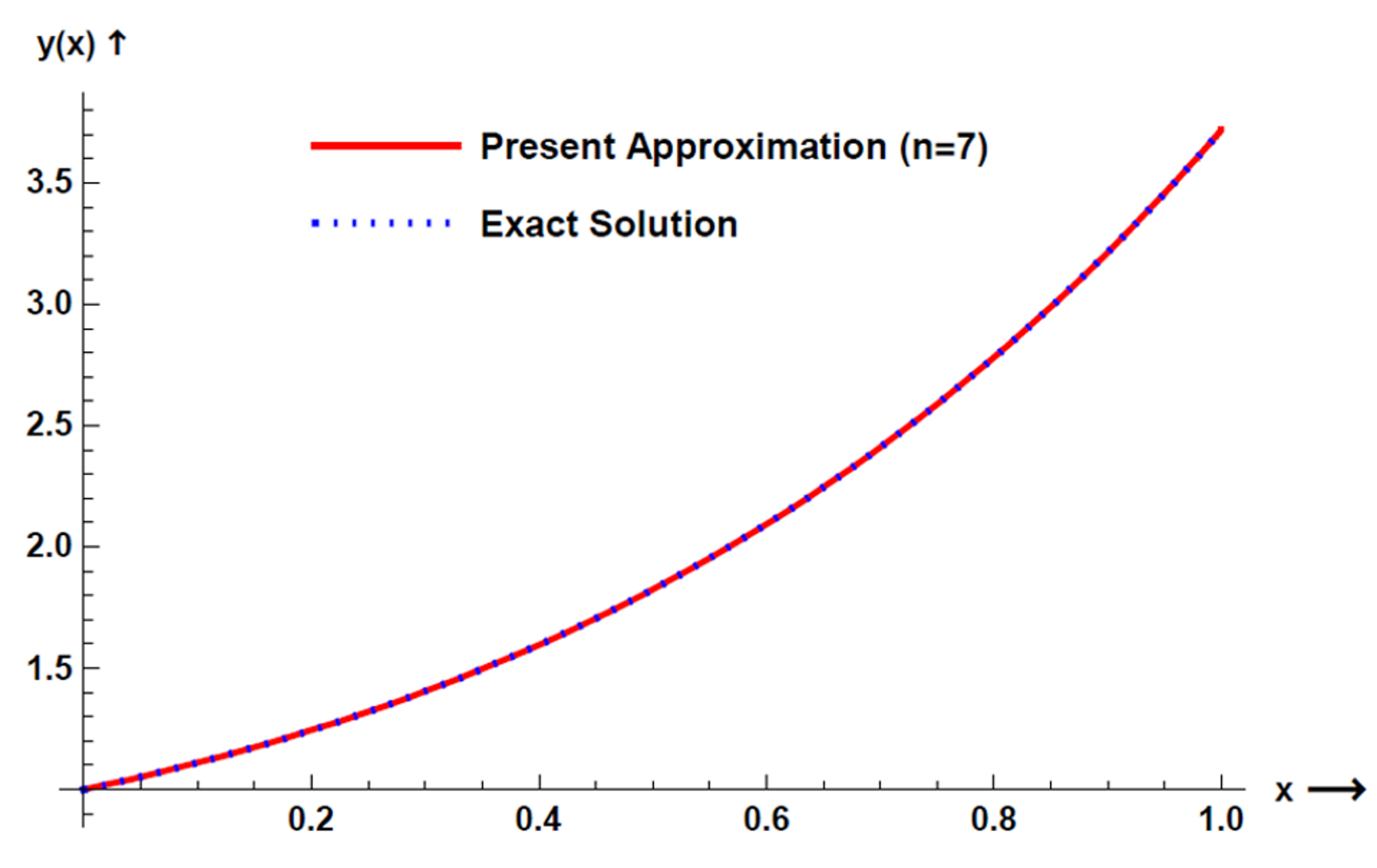} 
			\includegraphics[width=0.49\linewidth]{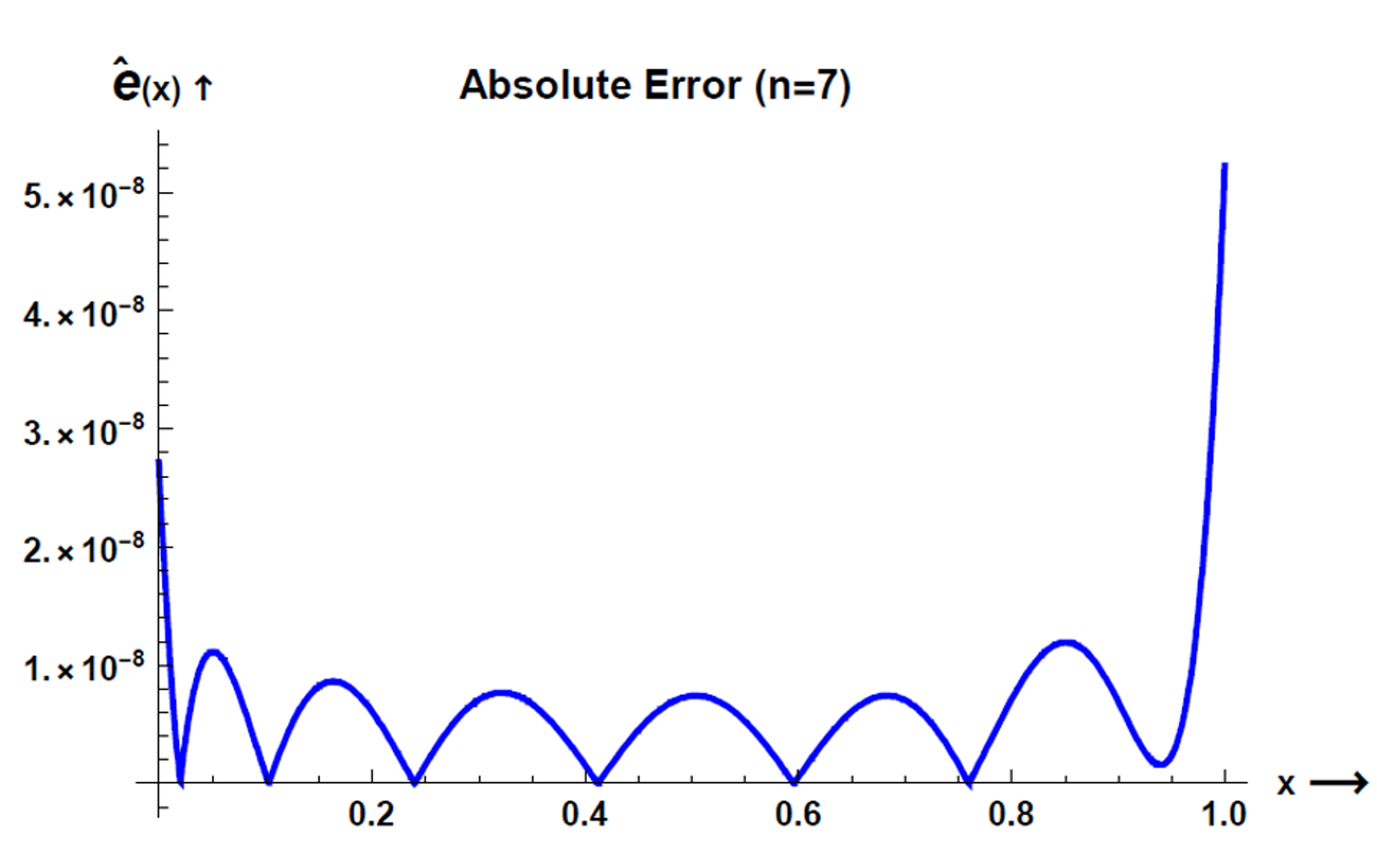}
			(a) \hspace{0.45\linewidth} (b)
			\caption{(a) Comparison of exact and present solution for example \ref{example 1} for $n=7$. (b) Absolute error between exact and approximate solutions of  \ref{example 1} for $n=7$.}
			\label{fig1}
		\end{figure}

	A comparison of approximation (\ref{y(x) for Ex.1}) with exact solution of problem  (\ref{example 1}) has been shown in figure \ref{fig1}. Maximum magnitude of the error between the two solutions is of order $10^{-8}$ for $n=7$.
	
	\begin{example} \label{example 2 }
		Let us consider the integro-differential equation (ref)
		
		\begin{equation} \label{Ex.2}
		\begin{split}
		(1+x^2)\frac{d^2y}{dx^2} + y + \cos x \int_{0}^{x}(x-t)^2\frac{dy}{dt}dt = r(x) \\
		y(0) = 0, \left(\frac{dy}{dx}\right)_{x=0} = 1; \:\: x \in [0,1]
		\end{split}
		\end{equation}
		
		For $r(x) = 2 (x - \sin x) \cos x - x^2 \sin x $, the problem (\ref{Ex.2}) admits the exact solution $y(x) =\sin x$.
	\end{example}

	As discussed in section (\ref{Case 2}), substituting 
	\begin{equation} \label{y" for ex.2}
	\frac{d^2y}{dx^2} = C^T \phi(x),
	\end{equation}
	
	\noindent into equation (\ref{Ex.2}) we get$-$ 
	\begin{equation}\label{simplified LHS ex.2}
	(1+x^2) C^T\phi(x) + C^T \Theta^2 \phi(x) +\cos x \, C^T \Theta^4 \phi(x) = r(x) - x - \frac{1}{3} x^3 \cos x
	\end{equation}
	
	Taking the approximations $(1+x^2)\phi(x) = A \phi(x), \, \cos x \phi(x) = B \phi(x), \, r(x) - x - \frac{1}{3} x^3 \cos x = R^T \phi(x)$, equation (\ref{simplified LHS ex.2}) can be simplified as:
	\begin{equation} \label{simplified expression for ex.2}
	C^T \left( A +  \Theta^2 + B\,\Theta^4  \right) \phi(x) = R^T \phi(x)
	\end{equation}
	\noindent where, $A$ and $B$ are matrices of order $n$, and $R$ is a column vector. For illustration, matrix $A$ can be calculated for $n=5$ as:
	\begin{equation*}
	A = \left[ {\begin{array}{*{20}{c}}
		{\frac{4}{3}}&{\frac{1}{{2\sqrt 3 }}}&{\frac{1}{{6\sqrt 5 }}}&0&0&0\\
		{\frac{1}{{2\sqrt 3 }}}&{\frac{7}{5}}&{\frac{1}{{\sqrt {15} }}}&\frac{1}{10}\sqrt {\frac{3}{7}}&0&0\\
		{\frac{1}{{6\sqrt 5 }}}&{\frac{1}{{\sqrt {15} }}}&{\frac{{29}}{{21}}}&{\frac{3}{{2\sqrt {35} }}}&{\frac{1}{{7\sqrt 5 }}}&0\\
		0&\frac{1}{10}\sqrt {\frac{3}{7}}&{\frac{3}{{2\sqrt {35} }}}&{\frac{{62}}{{45}}}&{\frac{2}{{3\sqrt 7 }}}&\frac{5}{9\sqrt{77}}\\
		{ - 30}&{ - \frac{{45\sqrt 3 }}{2}}&{ - \frac{{873}}{{14\sqrt 5 }}}&{ - \frac{{103}}{{3\sqrt 7 }}}&{ - \frac{{19}}{7}}&-\frac{5}{3\sqrt{11}}\\
		{\frac{{117\sqrt {11} }}{2}}&{43\sqrt {33} }&{114\sqrt {\frac{{11}}{5}} }&{\frac{1049}{18}\sqrt {\frac{11}{7}}}&{\frac{{17\sqrt {11} }}{3}}&{\frac{{29}}{9}}
	\end{array}} \right]
	\end{equation*}
	
\begin{figure}[h]
		\centering
		\includegraphics[width=0.49\linewidth]{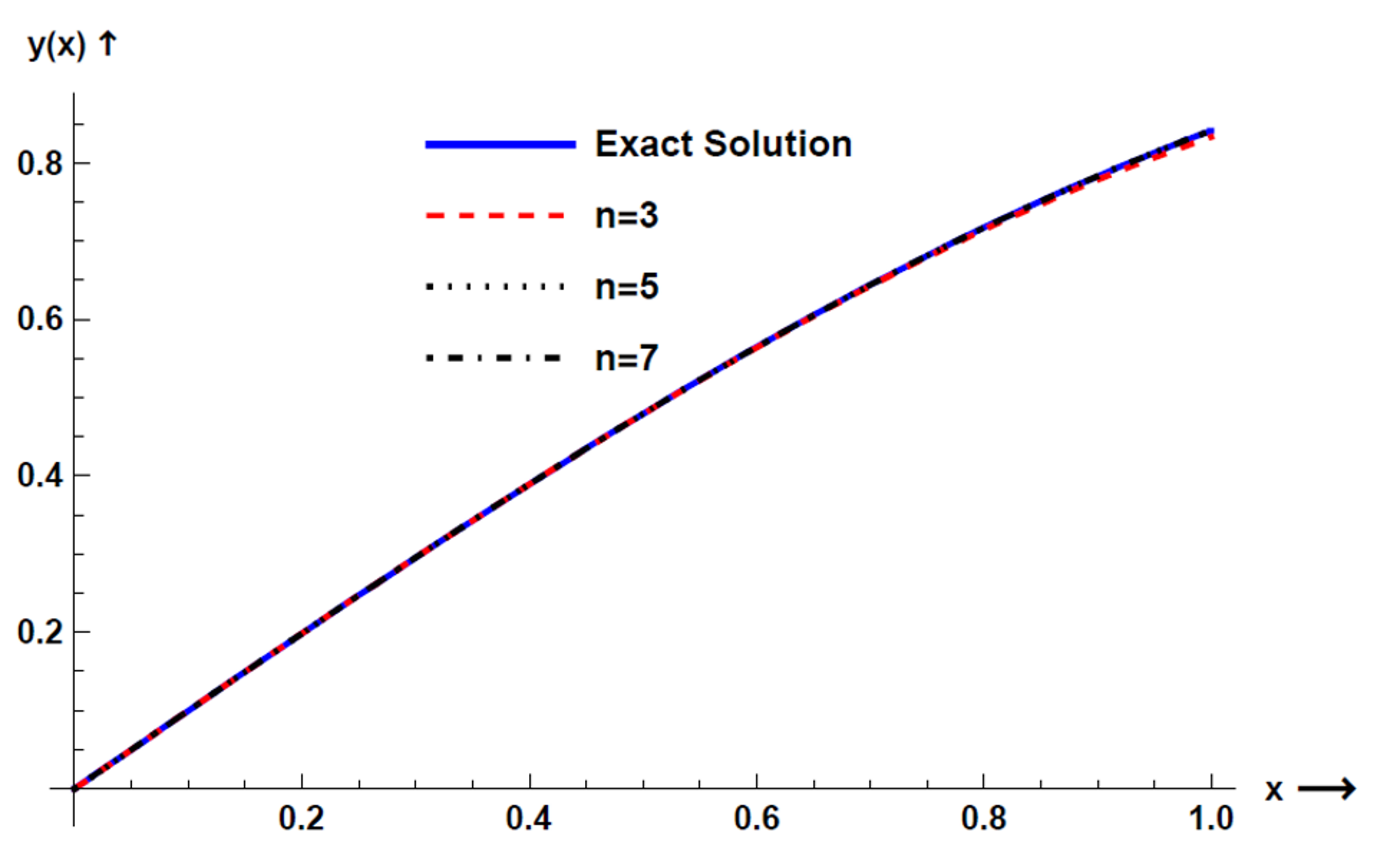} 
		\includegraphics[width=0.49\linewidth]{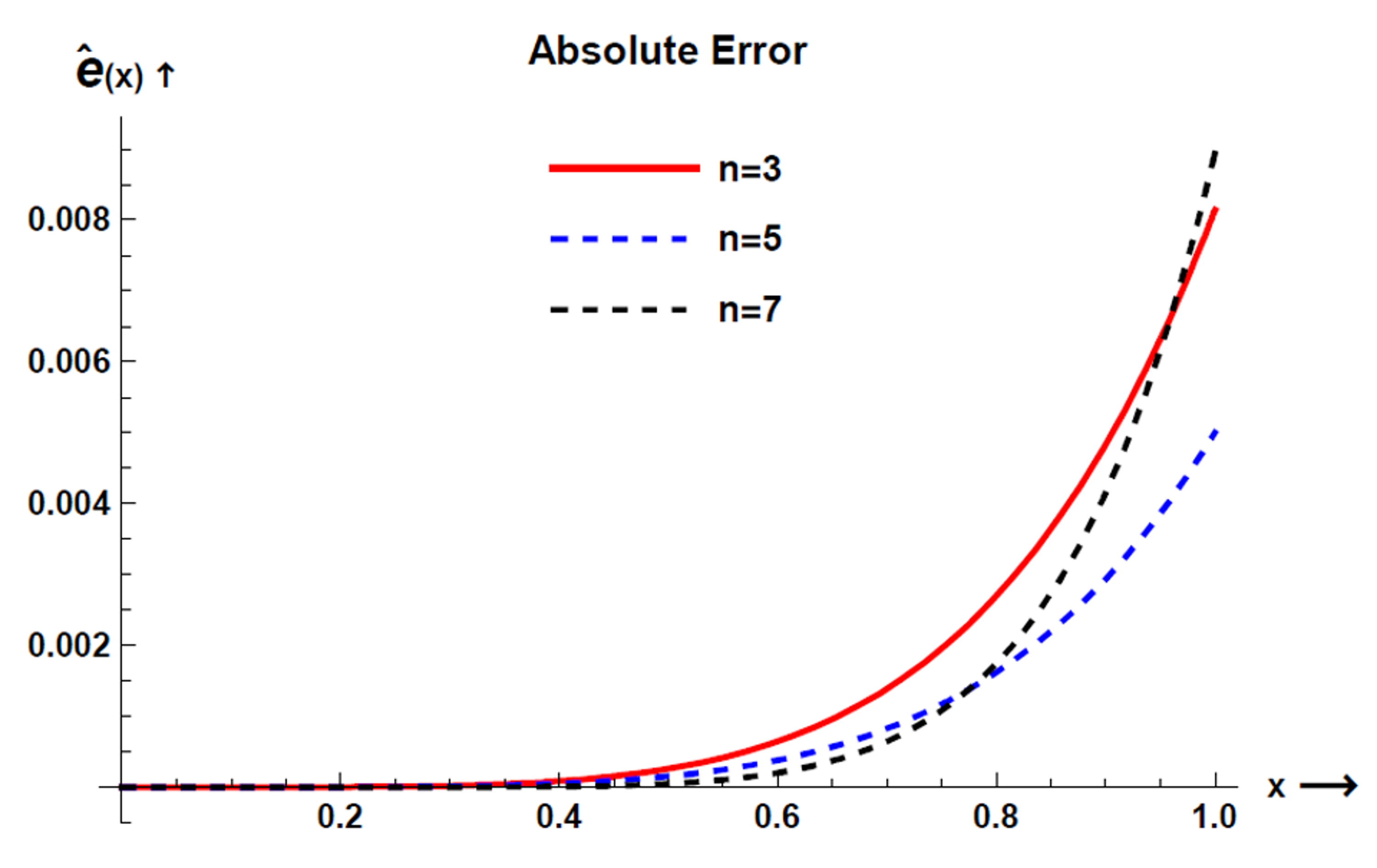}
		(a) \hspace{0.45\linewidth} (b)
		\caption{(a) Comparison of exact and present solution for example \ref{example 2 } for $n=3,5,7$. (b) Absolute error between exact and present solutions of  \ref{example 2 } for $n=3,5,7$.}
		\label{fig2}
	\end{figure}

	Thus, vector $C^T$ of unknown coefficients is obtained form equation (\ref{simplified expression for ex.2}) and, thereby, an approximation for $y(x)$ is obtained from equation (\ref{y" for ex.2}) as:  
	\begin{equation}
	y(x)_{n=3} \approx 0.99999998 x - 0.16667268 x^3
	\end{equation}
	\begin{equation} 
	y(x)_{n=5} \approx 0.99999998 x - 0.16667265 x^3 + 0.00003539 x^4 + 0.013110407 x^5
	\end{equation}
	\begin{equation} 
	\begin{split}
	y(x)_{n=7} \approx  0.99999998 x - 0.16667265 x^3 + 0.00003539 x^4 + 0.0083710487 x^5\\ + 
	0.0023319 x^6 - 0.011538172 x^7
	\end{split}
	\end{equation}
	Present approximations to solution of example (\ref{example 2 }) for $n = 3,5,7$ have been compared with the exact solution in figure \ref{fig2}.
	
	\section{Application to Ecology}
	In this section, we will apply this method to find solution of population problem (for females) \cite{Kot2001} 
	\begin{equation} \label{Ex.3}
	B'(t) = g(t) + \int_{0}^{t}K(t,\eta)B(\eta)d\eta 
	\end{equation}
	\noindent where, $B'(t)=\frac{dB}{dt}$, \\
	$K(t,\eta) = k(t-\eta)$ : net maternity function of females class age $\eta$ at time $t$.\\
	$g(t)$ : contribution of birth due to female already present at time $t$.\\
	$B(t)$ : the number of female births.
	
	Let the number of female births be given by $g(t) =\frac{1}{4} \left(6(1+t)- 7 e^{\frac{1}{2}t} - 4 \sin t\right)$, $K(t,\eta) = t- \eta$ and $t \in [0, 1]$. Then, the model (\ref{Ex.3}) takes the form:
	
	\begin{equation} \label{Ex.3 real}
	B'(t) -\int_{0}^{t}(t-\eta)B(\eta)d\eta  = e^t - \sin t; \: B(0)=1.
	\end{equation}
	This model admits the exact solution $B(t) = \frac{1}{2}\left(e^{\frac{1}{2}t} - \sin t + \cos t\right)$.
	
	As discussed in section \ref{section : Solution of IDE} and earlier examples, taking $B'(t)=C^T \phi(t)$, and $\frac{1}{2} t^2 +\frac{1}{4} \left(6(1+t)- 7 e^{\frac{1}{2}t} - 4 \sin t\right) =R^T \phi (t)$ for $n=5$, an approximate solution to equation (\ref{Ex.3 real}) is obtained as :
	
	\begin{equation}\label{y(x) for Ex 2}
	B(t) \approx 1 - 0.25002 t - 0.18730 t^2 + 0.09293 t^3 + 0.02375 t^4 - 
	0.00561 t^5
	\end{equation}
	
	The approximation (\ref{y(x) for Ex 2}) is compared with exact solution of population problem (\ref{Ex.3 real}) in figure \ref{fig2} for $n=5, 7$. The maximum magnitude of approximation errors for $n=5, 7$ are of order $10^{-5}$ and $10^{-7}$ respectively.
	
	\begin{figure}[h]
		\centering
		\includegraphics[width=0.49\linewidth]{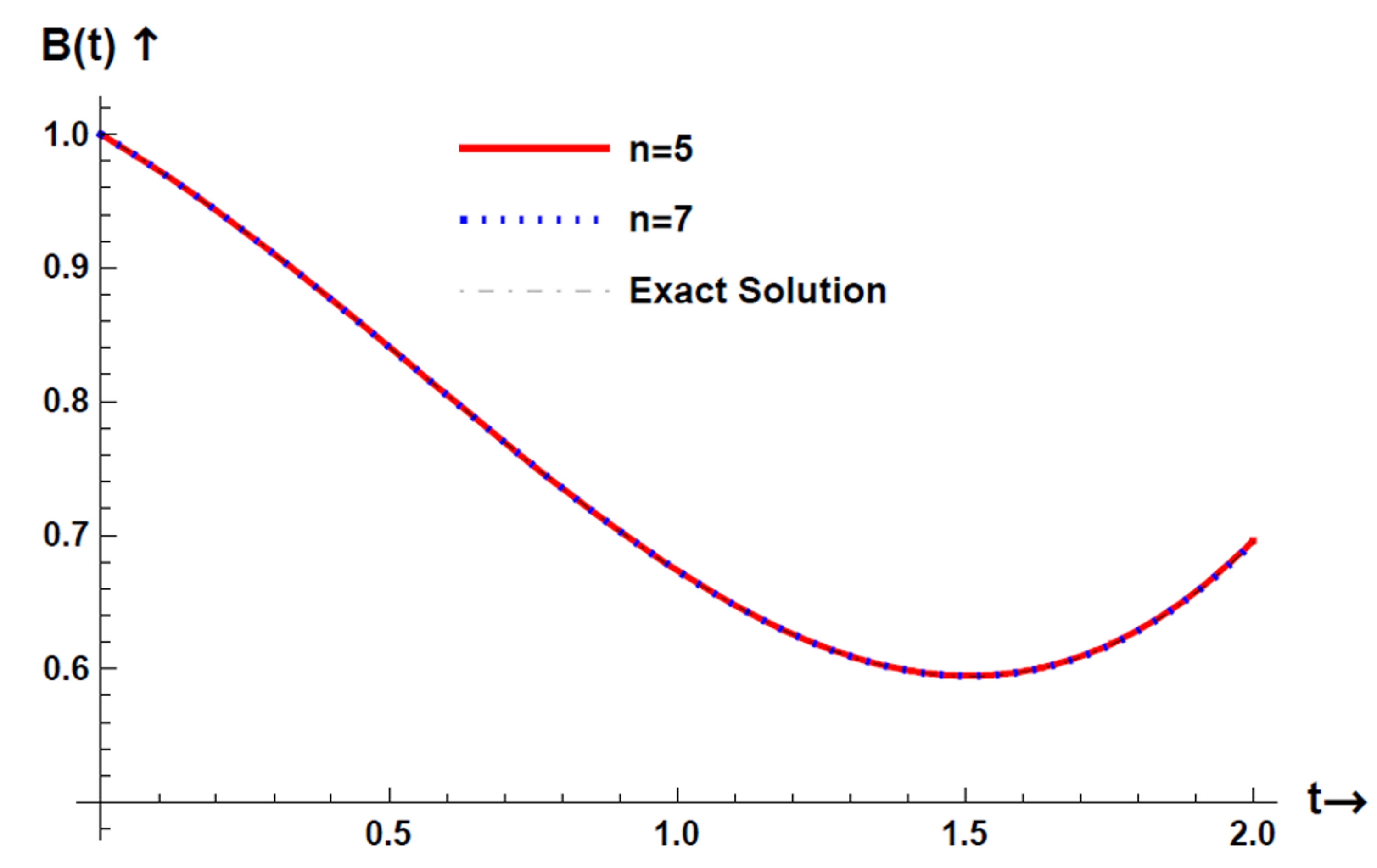} 
		\includegraphics[width=0.49\linewidth]{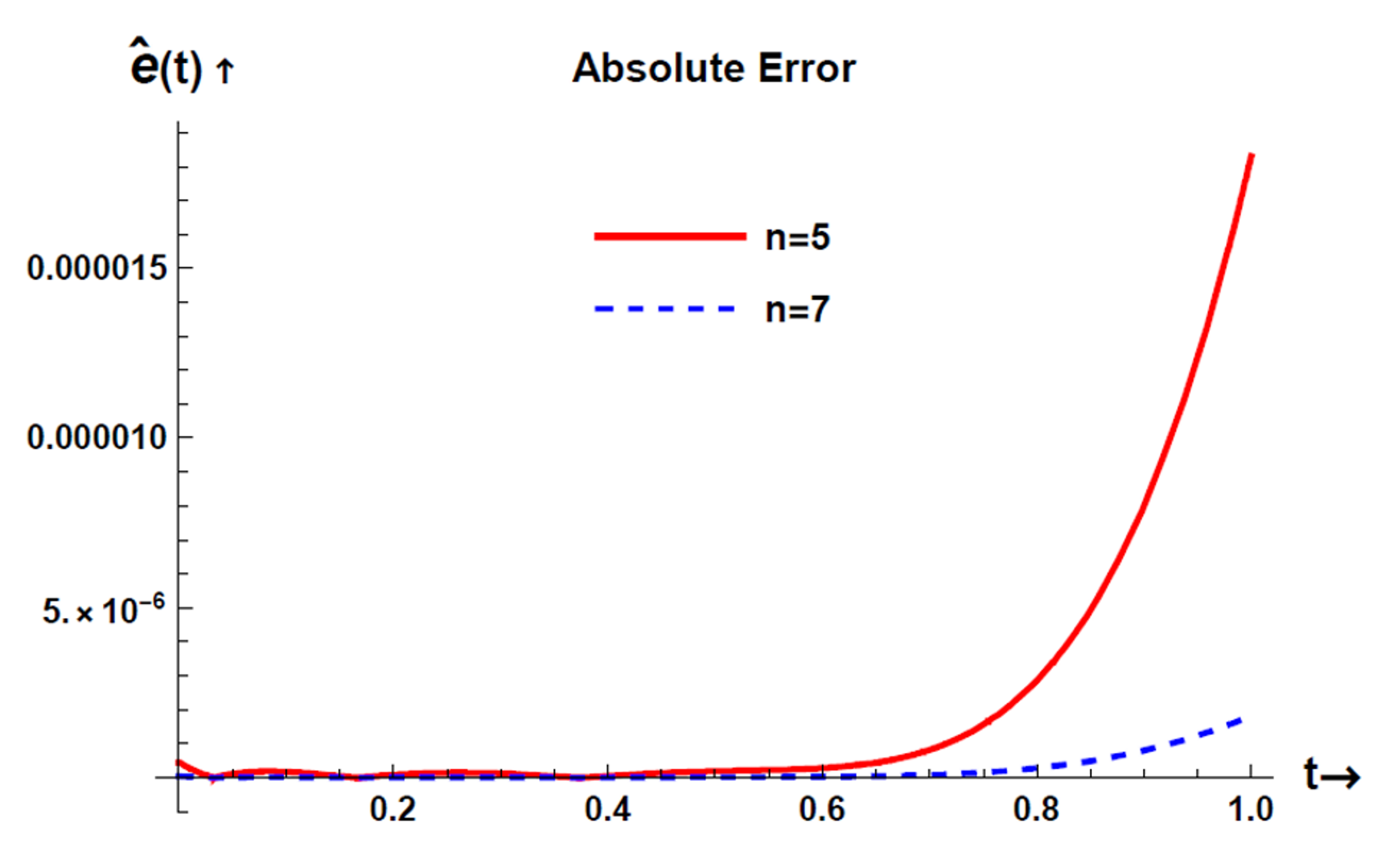}
		(a) \hspace{0.45\linewidth} (b)
		\caption{(a) Comparison of present approximation and exact solutions $(b)$ Absolute error $-$ between present approximation and exact solutions to population problem (\ref{Ex.3 real}) for $n=5,7$.}
		\label{fig3}
	\end{figure}
	
	\section{Conclusion}
	In this work, a new method has been applied to find approximate solution of linear integro-differential equations with help of Bernoulli polynomials. A set of $n$ orthonormal polynomials derived form Bernoulli polynomials of degree $n$ on $[0,1]$ has been used to form an operational matrix of integration. These new family of polynomials together with the operational matrix were applied to convert derivatives and integrals of dependent variable into an approximating polynomial form, thereby, converting an integro-differential equation into a set of algebraic equations with unknown coefficients, which are easily obtained with the help of operational matrix. Finally, an approximate solution is  obtained in form of a polynomial of degree $n$. Three integro-differential equations have been solved which includes one problem with constant coefficient, one with variable coefficient of general form and a population problem from earlier established literature. The problems have been solved for different values of $n$, numerical results have been compared with available exact solutions, and error of approximation have been plotted. It was concluded that most of the problems can be approximated by using only first few modified (orthonormal) polynomials with very small error. Some outcomes of this method can be summarized as follows.
	\begin{itemize}
		\item error is small and can be reduced by taking higher degree approximation. 
		\item the method is very fast for integro-differential equations with constant coefficients.
		\item solution is obtained in form of a polynomial, which can be easily carried forward for various further applications. 
		\item error can be minimized up to required accuracy because error decreases quickly with increase of $n\, -\, $the degree of Bernoulli polynomials.
		\item method can be programmed for various numerical applications.
	\end{itemize}

\bibliographystyle{ieeetr}
\bibliography{bibliography}

\begin{thebibliography}{10}

\bibitem{McKendrick1926}
A.~G. McKendrick, ``Applications of mathematics to medical problems,'' {\em
  Proc. Edinb. Math. Soc.}, vol.~44, pp.~98--130, 1926.

\bibitem{Schaffer1986}
M.~Kot and W.~M. Schaffer, ``{Discrete-time growth–dispersal models},'' {\em
  Math. Biosci.}, vol.~80, no.~1, pp.~283--326, 1986.

\bibitem{Kendall1965}
D.~G. Kendall, ``{Mathematical models of the spread of infection},'' {\em
  Mathematics and Computer Science in Biology and Medicine}, vol.~80, no.~1,
  pp.~213--225, 1965.

\bibitem{Baym1990}
G.~Baym, {\em Lectures on Quantum Mechanics}.
\newblock Redwood City, California: Addison–Wesley, 1990.

\bibitem{Vangipuram1995}
V.~Lakshmikantham and M.~R.~M. Rao, {\em Theory of Integro-Differential
  Equations}.
\newblock Switzerland: Gordon and Breach Science Publisher, 1978.

\bibitem{Petrin1997}
A.~B. Petrin, ``{Integro-differential equation method in the hydrodynamics of
  an incompressible viscous fluid},'' {\em Journal of Experimental and
  Theoretical Physics}, vol.~84, no.~4, pp.~724--727, 1997.

\bibitem{Kot2001}
M.~Kot, {\em Elements of Mathematical Ecology}.
\newblock Cambridge University Press, 2001.

\bibitem{Darania2007}
P.~Darania and A.~Ebadian, ``{A method for the numerical solution of the
  integro-differential equations},'' {\em Applied Mathematics and Computation},
  vol.~188, no.~1, pp.~657--668, 2007.

\bibitem{Abdelaziz2010}
A.~Mennouni and S.~Guedjiba, ``{A note on solving integro-differential equation
  with Cauchy kernel},'' {\em Mathematical and Computer Modelling}, vol.~52,
  no.~9-10, pp.~1634--1638, 2010.

\bibitem{Yuzbasi2011}
u.~Y{\"{u}}zbaşi, ``{A numerical approach for solving the high-order linear
  singular differential-difference equations},'' {\em Computers and Mathematics
  with Applications}, vol.~62, no.~5, pp.~2289--2303, 2011.

\bibitem{Sweilam2013}
N.~H. Sweilam, M.~M. Khader, and W.~Y. Kota, ``{Numerical and analytical study
  for fourth-order integro-differential equations using a pseudospectral
  method},'' {\em Mathematical Problems in Engineering}, vol.~2013, no.~Article
  ID 434753, pp.~1--7, 2013.

\bibitem{Dehghan2012}
M.~Dehghan and R.~Salehi, ``{The numerical solution of the non-linear
  integro-differential equations based on the meshless method},'' {\em Journal
  of Computational and Applied Mathematics}, vol.~236, no.~9, pp.~2367--2377,
  2012.

\bibitem{Maleknejad2012}
K.~Maleknejad, B.~Basirat, and E.~Hashemizadeh, ``{A Bernstein operational
  matrix approach for solving a system of high order linear Volterra-Fredholm
  integro-differential equations},'' {\em Mathematical and Computer Modelling},
  vol.~55, no.~3-4, pp.~1363--1372, 2012.

\bibitem{Yuzbasi2012}
u.~Y{\"{u}}zbaşi, N.~Şahin, and A.~Yildirim, ``{A collocation approach for
  solving high-order linear Fredholm-Volterra integro-differential
  equations},'' {\em Mathematical and Computer Modelling}, vol.~55, no.~3-4,
  pp.~547--563, 2012.

\bibitem{Yuzbasi2013}
u.~Y{\"{u}}zbaşI and M.~Sezer, ``{An improved Bessel collocation method with a
  residual error function to solve a class of Lane-Emden differential
  equations},'' {\em Mathematical and Computer Modelling}, vol.~57, no.~5-6,
  pp.~1298--1311, 2013.

\bibitem{Mirzaee2015}
F.~Mirzaee and S.~Bimesl, ``{Numerical solutions of systems of high-order
  Fredholm integro-differential equations using Euler polynomials},'' {\em
  Applied Mathematical Modelling}, vol.~39, no.~22, pp.~6767--6779, 2015.

\bibitem{Singh2017b}
S.~Singh, V.~K. Patel, V.~K. Singh, and E.~Tohidi, ``{Numerical solution of
  nonlinear weakly singular partial integro-differential equation via
  operational matrices},'' {\em Applied Mathematics and Computation}, vol.~298,
  pp.~310--321, 2017.

\bibitem{Jafarzadeh2018}
Y.~Jafarzadeh and B.~Keramati, ``{Numerical method for a system of
  integro-differential equations and convergence analysis by Taylor
  collocation},'' {\em Ain Shams Engineering Journal}, vol.~9, no.~4,
  pp.~1433--1438, 2018.

\bibitem{Katsikadelis2019}
J.~T. Katsikadelis, ``{Numerical solution of integrodifferential equations with
  convolution integrals},'' {\em Archive of Applied Mechanics}, vol.~89,
  no.~10, pp.~2019--2032, 2019.

\bibitem{Singh2020}
Y.~Singh, V.~Gill, J.~Singh, D.~Kumar, and K.~S. Nisar, ``{On the Volterra-Type
  Fractional Integro-Di ff erential Equations Pertaining to Special
  Functions},'' {\em Fractal and fractional}, vol.~4, no.~33, pp.~1--12, 2020.

\bibitem{Cheon2003}
G.~S. Cheon, ``{A note on the Bernoulli and Euler polynomials},'' {\em Applied
  Mathematics Letters}, vol.~16, no.~3, pp.~365--368, 2003.

\bibitem{Maleknejad2007}
K.~Maleknejad, S.~Sohrabi, and Y.~Rostami, ``{Numerical solution of nonlinear
  Volterra integral equations of the second kind by using Chebyshev
  polynomials},'' {\em Applied Mathematics and Computation}, vol.~188, no.~1,
  pp.~123--128, 2007.

\bibitem{Nemati2015}
S.~Nemati, ``{Numerical solution of Volterra-Fredholm integral equations using
  Legendre collocation method},'' {\em Journal of Computational and Applied
  Mathematics}, 2015.

\bibitem{Rahman2012}
M.~A. Rahman, M.~S. Islam, and M.~M. Alam, ``{Numerical Solutions of Volterra
  Integral Equations Using Laguerre Polynomials},'' {\em Journal of Scientific
  Research}, vol.~4, no.~2, pp.~357--364, 2012.

\bibitem{Yousefi2006}
S.~A. Yousefi, ``{Numerical solution of Abel's integral equation by using
  Legendre wavelets},'' {\em Applied Mathematics and Computation}, vol.~175,
  no.~1, pp.~575--580, 2006.

\bibitem{Sahu2019}
P.~K. Sahu and B.~Mallick, ``{Approximate Solution of Fractional Order
  Lane–Emden Type Differential Equation by Orthonormal Bernoulli's
  Polynomials},'' {\em International Journal of Applied and Computational
  Mathematics}, vol.~5, no.~89, 2019.

\bibitem{Tohidi2013}
E.~Tohidi, A.~H. Bhrawy, and K.~Erfani, ``{A collocation method based on
  Bernoulli operational matrix for numerical solution of generalized pantograph
  equation},'' {\em Applied Mathematical Modelling}, vol.~37, no.~6,
  pp.~4283--4294, 2013.

\bibitem{Tohidi2013a}
E.~Tohidi and A.~Kili{\c{c}}man, ``{A collocation method based on the bernoulli
  operational matrix for solving nonlinear BVPs which arise from the problems
  in calculus of variation},'' {\em Mathematical Problems in Engineering},
  vol.~2013, no.~Article ID 757206, pp.~1--9, 2013.

\bibitem{Mohsenyzadeh2016}
M.~Mohsenyzadeh, ``{Bernoulli operational Matrix method of linear Volterra
  integral equations},'' {\em Journal of Industrial Mathematics}, vol.~8,
  no.~3, pp.~201--207, 2016.

\bibitem{Costabile2006}
F.~A. Costabile and F.~Dell'Accio, ``{A new approach to Bernoulli
  polynomials},'' {\em Rendiconti di Matematica, Serie VII}, vol.~26,
  pp.~1--12, 2006.

\bibitem{Todorov1984}
P.~G. Todorov, ``{On the theory of the Bernoulli polynomials and numbers},''
  {\em Journal of Mathematical Analysis and Applications}, vol.~104, no.~2,
  pp.~309--350, 1984.

\bibitem{Kreyszig1978}
K.~E., {\em Introductory Functional Analysis with Applications}.
\newblock New York, USA: John Wiley and Sons Press, 1978.

\bibitem{Kurt2011}
B.~Kurt and Y.~Simsek, ``{Notes on generalization of the Bernoulli type
  polynomials},'' {\em Applied Mathematics and Computation}, vol.~218, no.~3,
  pp.~906--911, 2011.

\bibitem{Natalini2003}
P.~Natalini and A.~Bernardini, ``{A generalization of the Bernoulli
  polynomials},'' {\em Journal of Applied Mathematics}, vol.~3, no.~3,
  pp.~155--163, 2003.

\bibitem{Lu2011}
D.~Q. Lu, ``{Some properties of Bernoulli polynomials and their
  generalizations},'' {\em Applied Mathematics Letters}, vol.~24, no.~5,
  pp.~746--751, 2011.

\end{thebibliography}

\end{document}